\documentclass[fleqn,final]{zzz}

\usepackage{mathrsfs}
\usepackage{color}
\usepackage{tikz}
\usepackage{amsmath,amsthm,amsbsy,amssymb}
\usepackage{graphicx}
\usepackage{rotating}
\usepackage{fixmath}
\usepackage[pdftex]{hyperref}
\usepackage{soul}
\usepackage{enumitem}

\hypersetup{
 colorlinks = true,
 urlcolor = blue,
 linkcolor = red,
 citecolor = green,
}

\sethlcolor{black}

\makeatletter
\def\eqnarray{\stepcounter{equation}\let\@currentlabel=\theequation
\global\@eqnswtrue
\tabskip\@centering\let\\=\@eqncr
$$\halign to \displaywidth\bgroup\hfil\global\@eqcnt\z@
$\displaystyle\tabskip\z@{##}$&\global\@eqcnt\@ne
\hfil$\displaystyle{{}##{}}$\hfil
&\global\@eqcnt\tw@ $\displaystyle{##}$\hfil
\tabskip\@centering&\llap{##}\tabskip\z@\cr}

\def\endeqnarray{\@@eqncr\egroup
\global\advance\c@equation\m@ne$$\global\@ignoretrue}

\def\@yeqncr{\@ifnextchar [{\@xeqncr}{\@xeqncr[5pt]}}
\makeatother

\usepackage{scalerel,url}
\let\svus_
\catcode`_=\active %
\def_{\ifmmode\expandafter\svus\expandafter\bgroup\expandafter\lowerit
\else\svus\fi}
\def\lowerit#1{\ThisStyle{\raisebox{-2\LMpt}{$\SavedStyle#1$}}\egroup}
\makeatletter
 \AtBeginDocument{
 \check@mathfonts
 \fontdimen13\textfont2=3.5pt
 \fontdimen14\textfont2=3.5pt
 \fontdimen16\textfont2=2.5pt
 \fontdimen17\textfont2=2.5pt
 }
\makeatother

\newcommand{\bbC}{{\mathbb{C}}}
\newcommand{\bbD}{{\mathbb{D}}}

\newcommand{\bbR}{{\mathbb{R}}}

\newcommand{\fre}{{\frak{e}}}

\newcommand{\x}{{\mathbf{x}}}

\newcommand{\beq}{\begin{equation}}
\newcommand{\eeq}{\end{equation}}
\newcommand{\ba}{\begin{align}}
\newcommand{\ea}{\end{align}}

\newcommand{\norm}[1]{\lVert#1\rVert}
\newcommand{\overbar}{\overline}

\begin{document}

\renewcommand{\PaperNumber}{***}

\FirstPageHeading

\ShortArticleName{Commentary on Askey's Szeg\H{o} paper}

\ArticleName{Commentary on Askey's Szeg\H{o} paper}

\Author{Barry Simon$\,^{\ast}$
}

\AuthorNameForHeading{B.~Simon}
\Address{$^\ast$ IBM Professor of Mathematics and Theoretical Physics, Emeritus,
California Institute of Technology, Pasadena, California  91125, USA
\URLaddressD{
\href{https://en.wikipedia.org/wiki/Barry_Simon}
{https://en.wikipedia.org/wiki/Barry\_Simon}
}   


} 
\EmailD{bsimonca@gmail.com} 




\ArticleDates{Received \today~in final form ????; Published online ????}



\vspace{0.5cm}

\noindent One of the celebrated gems of Dick Askey is his editing of the complete works of G\'abor Szeg\H{o}. Dick's commentary on the individual papers are often brilliant expositions of the subject of that paper of Szeg\H{o}, sometimes the general context rather than the content of the paper itself.\\[-0.2cm]

\noindent My contribution to Dick's Selecta is in this spirit. Since I am writing about Dick's take on Szeg\H{o}'s work, this is about my thinking on that same subject.\\[-0.2cm]

\noindent Our views are different. The modern theory of orthogonal polynomials is usually broken into two parts, often called the analytic theory and the algebraic theory. The former deals with the process of going from a general measure on $\bbC$ with finite moments (most often measures on $\bbR$ or $\partial\bbD$ called orthogonal polynomials on the real line or unit circle, aka OPRL or OPUC) to its orthogonal polynomials which uses tools of mathematical analysis. The latter deals with specific special cases going back to the so called classical families of Jacobi, Laguerre and Hermite. Here the techniques are typically algebraic although occasionally analytic.\\[-0.2cm]

\noindent Szeg\H{o} was the acknowledged master in his generation (and I think all generations) of both sides (see below). Dick was the acknowledged master in his generation of the algebraic side while my work has been on the analytic side, especially using methods of spectral theory \cite{OPUC1, OPUC2}.\\[-0.2cm]

\noindent Askey's discussion of Szeg\H{o}'s work begins with Szeg\H{o}'s first paper which solved a problem of P\'{o}lya. This result, called the Szeg\H{o} limit theorem, and its extensions, is one that Szeg\H{o} returned to often. To many, including me, the result is his most important -- it should be mentioned that Askey's paper states that Szeg\H{o}'s most important work is his problem book with P\'{o}lya \cite{PS} (P\'{o}lya also expressed this opinion). Askey also quotes the opinion of the physicist, Barry McCoy, that Szeg\H{o}'s most important work was the strong Szeg\H{o} limit theorem (I will shortly explain what the limit theorems are and will later explain why I disagree with Askey and McCoy).\\[-0.2cm]

\noindent P\'{o}lya's conjecture involves Toeplitz determinants. Given a positive $L^1$ function,  $f$, on the unit circle, the degree $n$ Toeplitz matrix, $T_n(f)$, is the $n\times n$ matrix whose $kj$ matrix element is the $k-j$ Fourier coefficient of $f$, $\hat{f}(k-j)$, where

\begin{equation}\label{1}
  \hat{f}(\ell) = \int_0^{2\pi}  f(\theta) \exp(-i\ell\theta) \frac{d\theta}{2\pi}.
\end{equation}

\noindent (note in terms of the natural variable, $z= \exp(i\theta)$, this is just the (conjugate conjugates of the) moments of the measure, $f\tfrac{d\theta}{2\pi}$). The determinant of $T_{n+1}(f)$ is denoted $D_n(f)$ (note the $n+1$; $D_n$ is thus an $n\times$ determinant, so $D_0=f(0)=1$ when $f$ is normalized). P\'{o}lya conjectured that $D_n^{1/n} $ had a limit. What Szeg\H{o} proved is not only that the limit exists but its value, namely that

\begin{equation}\label{2}
\lim_{n\to\infty} \frac{\log(D_n(f))}{n}  =   \int_0^{2\pi} \log(f(\theta)) \frac{d\theta}{2\pi}.
\end{equation}

\noindent This is the Szeg\H{o} limit theorem. The usual historical reference for this result is \cite{Sz2} although Askey quotes \cite{Sz1}, which is an announcement.\\[-0.2cm]

\noindent The modern theory considers not positive functions but probability measures, $\mu$, on the unit circle of the form

\begin{equation}\label{3}
  d\mu =  f(\theta) \frac{d\theta}{2\pi} + d\mu_s,
\end{equation}

\noindent where the latter is singular with respect to Lebesgue measure on the unit circle.  Remarkably it is known that (2) always holds although both sides may be $-\infty$. Note that the limit is independent of the singular part of the measure.\\[-0.2cm]

\noindent As mentioned, Szeg\H{o} returned to this subject often. I want to focus on two times. One was in 1952 when he found (\cite{Sz5}) the second term in an asymptotic series for $\log(D_n) $. To state the result in general form one considers only measures for which $\log(f)$ is an $L^1$ function on the circle (which is equivalent to $\lim_n D_n^{1/n}>0$), in which case one uses $L_n$ for its Fourier coefficients. Then it is known that whenever

\begin{equation}\label{4}
 \sum_{n}  n |L_n|^2<\infty
\end{equation}

\noindent one has that

\begin{equation}\label{5}
  \lim_{n\to\infty} [\log(D_n) - n L_0] =L_0 + \sum_{n>0} n |L_n|^2.
\end{equation}

\noindent This is called the strong Szeg\H{o} theorem although Szeg\H{o} only proved it under some restrictions on the measure (including that $\mu_s = 0$). My book \cite{OPUC1} has several proofs of the Szeg\H{o} limit theorem in Chapter 2 and of the strong Szeg\H{o} limit theorem in Chapter 6.\\[-0.2cm]

\noindent Remarkably, it was almost 40 years between finding the first and second term of the series.  It was not because Szeg\H{o} searched for this second term hard for many years and there was a breakthrough. Indeed, one has the impression that until the late 1940's Szeg\H{o} had not considered the problem.  Instead the paper sprang
from a conjecture of the Yale chemist Lars Onsager.  While Onsager is best known among mathematicians and theoretical physicists for his exact solution of the Ising model, he got the Nobel Prize in chemistry for unrelated work in thermodynamics.  In connection with the solution of the Ising model, Onsager came across some explicit Toeplitz determinants and needed some conjectures on the second term in the asymptotics of this special case.  He asked the celebrated Yale mathematician, Kakutani, who knew that the person to bring it to was Szeg\H{o}.\\[-0.2cm]

\noindent Undoubtedly, the most important outgrowth of Szeg\H{o}'s limit theorem is  the theory of OPUC.  While Jacobi is a key player in a general theory of OPRL, that subject has several fathers, but OPUC is the invention of Szeg\H{o} alone, primarily in a pair of his papers \cite{Sz3}. Notice that they are not titled in terms of orthogonal polynomials but in terms of Toeplitzschen Formen. Given a positive measure on the unit circle (this was before general measures were widely considered and Szeg\H{o} considered only purely absolutely continuous measures), one gets monic orthogonal polynomials, $\{\Phi_n(z)\}_{n=0}^\infty$ by applying Gram--Schmidt to $\{z^n\}_{n=0}^\infty$.  What Szeg\H{o} found was that if $D_n$ are the Toeplitz determinants associated to the Fourier coefficients of the measure, then, for $n\ge 1$

\begin{equation}\label{6}
    \norm{\Phi_n}^2=\frac{D_n}{D_{n-1}}
\end{equation}

\noindent so that, since $D_0=1$, one has that

\begin{equation}\label{6A}
    D_n=\prod_{j=1}^{n}  \norm{\Phi_j}^2
\end{equation}

\noindent Thus, his limit theorem says something about asymptotics of norms.  I note that the \eqref{6A} is sometimes called Heine's formula because Heine proved the analog for OPRL.  You can find proofs and further discussion in \cite[Sections 1.2.7, 1.5]{OPUC1}. \\[-0.2cm]

\noindent In 1939, when he gave the Colloquium Lectures of the AMS, which appeared as the first edition of \cite{Sz4}, Szeg\H{o} revisited the theory of OPUC (which got little attention until it became a Russian sport) and, in particular, wrote down the recursion relation for the OPUC, namely

\begin{equation}\label{7}
  \Phi_{n+1}(z) = z\Phi_n(z) - \overbar{\alpha}_n\Phi_{n}^*(z); \quad \Phi_n^*(z) \equiv z^n\overline{\Phi_n(1/\bar{z})},
\end{equation}

\noindent for a suitable sequence of constants, $\{\alpha_n\}_{n=0}^\infty$ each in $\bbD$.  These parameters, written down first by Szeg\H{o}, are called Verblunsky coefficients (see below).  They are the OPUC analog of the Jacobi parameters of OPRL.  A derivation of \eqref{7} can be found, for example, in \cite[Section 1.5]{OPUC1}.\\[-0.2cm]

\noindent The derivation of \eqref{7} depends on noting that $\Phi_{n+1}(z)-z\Phi_n(z)$ is orthogonal to $\Phi_{n+1}$ which is the start of the proof of \eqref{7}, and which implies that $\norm{\Phi_{n+1}}^2+|\alpha_n|^2\norm{\Phi_n}^2 = \norm{\Phi_n}^2$ which, in turn, implies that

\begin{equation}\label{8}
  1-|\alpha_n|^2 = \frac{\norm{\Phi_{n+1}}^2}{\norm{\Phi_n}^2}
\end{equation}

By the same argument that led to \eqref{6A}, this implies

\begin{equation}\label{8A}
    \norm{\Phi_n}^2 = \prod_{j=0}^{n-1} (1-|\alpha_j|^2)
\end{equation}

\noindent The right side is clearly monotone decreasing in $n$ so both sides have limits and we see that

\begin{equation}\label{8B}
    \lim_{n\to\infty} \norm{\Phi_n}^2 = \prod_{j=0}^{\infty} (1-|\alpha_j|^2)
\end{equation}

\noindent By \eqref{6}, we see that $\lim_{n\to\infty}\tfrac{D_n}{D_{n-1}}$ exists and is equal to either side of \eqref{8B}.  \eqref{8A} implies the norm is less than one and that it (and so, by \eqref{6}, the ratio) is monotone decreasing.   From this, it is not hard to see that $\lim_{n\to\infty}D_n^{1/n}$ exists and is equal to $\lim_{n\to\infty}\tfrac{D_n}{D_{n-1}}$ and to either side of \eqref{8B}.  We thus have a sketch of a solution of P\'{o}lya's problem that the limit exists but not the full Szeg\H{o} limit theorem that expresses the limit in terms of an integral of $\log(f)$.

If we put together \eqref{8B}, the argument just given that the limit of $D_n^{1/n}$ is either side of \eqref{8B} and the Szeg\H{o} limit theorem in the form \eqref{2}, we get Verblunsky's form of the Szeg\H{o} limit theorem

\begin{equation}\label{9}
  \prod_{n=0}^\infty (1-|\alpha_n|^2) = \exp\left( \int_0^{2\pi}  \log(f(\theta)) \frac{d\theta}{2\pi}\right).
\end{equation}

\noindent Killip and I \cite{KiSi} found what is an OPRL analog of this expression.  We learned from the Russian literature that this formula was a form of Szeg\H{o}'s limit theorem.  I got interested in where this version of the theorem first appeared and I learned it was in an obscure paper of Verblunsky \cite{Verb} who didn't have the full Szeg\H{o} recursion formula but defined parameters equivalent to Szeg\H{o}'s $\alpha_n$.  Verblunsky's paper was little noted at the time.  While, in my opinion, the $\alpha_n$ should be called Szeg\H{o} coefficients, at the time I got interested in the subject (shortly after the turn of the century), they (or $-\alpha_n$ or $\pm\bar{\alpha}_n$) had three or four common names, one of which was Szeg\H{o} coefficients.  Seeking to recognize Verblunsky's work and realizing the advantage of starting in a new direction, I made a push to use Verblunsky coefficients and due to the force of my personality and the impact of my OPUC books \cite{OPUC1, OPUC2}, the name has stuck.\\[-0.2cm]

\noindent I promised to return to the issue of what is Szeg\H{o}'s greatest work.  No doubt, his problem book with P\'{o}lya is wonderful, arguably the best problem book ever.  But I share the opinion of many research mathematicians (even though I am proud of many of my own books) that one's original contributions count more than the most brilliant expositions of overall subjects.  This is not to say that there isn't original mathematics in P\'{o}lya--Szeg\H{o}.  As for McCoy's opinion about the strong Szeg\H{o} limit theorem, he is a theoretical physicist, best known for his work on Onsager's solution, so it is natural for him to most highly prize that work of Szeg\H{o} that is close to his interests.  As an expert on OPUC, I am, of course, not unbiased. That said, the leading term theorem makes the strong limit theorem possible and is central to large areas of analysis.\\[-0.2cm]

\noindent As the penultimate topic, I want to note that while I have restricted the discussion to Szeg\H{o}'s two limit theorems, mainly because that is what Askey's note \cite{Ask1} discusses, Szeg\H{o} has many other significant results.  Three of my favorites are his work on
\begin{enumerate}[label=(\roman*)]
\item Chebyshev polynomials and potential theory \cite{Sz1a} (the so-called Faber--Fekete--Szeg\H{o} theorem says that for any compact subset, $\fre\subset\bbC$, if $\norm{\cdot}_\fre$ is the sup norm on $\fre$ and $T_n^{(\fre)}$ the $n$th Chebyshev polynomial of $\fre$, then $\lim_{n\to\infty}\norm{T_n^{(\fre)}}_\fre^{1/n}=C(\fre)$, the logarithmic capacity of $\fre$);
\item his lovely theorem \cite{Sz2a} on power series with coefficients taking only finitely many values (that if $f(z)=\sum_{n=1}^\infty c_nz^n$ where $\{c_n\}_{n=1}^\infty$ takes only a finite set of values, then the analytic function defined on the disk by the power series either has a natural boundary on the entire unit disk or else is a rational function with poles only at all the $k$th roots of unity in which case the $c_n$ are eventually periodic of period $k$); and
\item  his foundational results on Hardy spaces (the existence \cite{Sz3} of the Szeg\H{o} reproducing kernel for $H^2$ which plays an important role in higher dimensional analysis and his proof \cite{Sz3a}, independently of Fatou, that if $f(\exp(i\theta))$ is the a.e.~boundary value of an $H^p$ function, then $\log(f(\exp(i\theta)))$ lies in $L^1$; this result provides the most elegant proof that the set of zeros of such functions has Lebesgue measure zero).
\end{enumerate}
\noindent  Notice that in these papers both $n$th root asymptotics and $\int \log(f)$ which were central to Szeg\H{o}'s first paper recur.\\[-0.2cm]

\noindent In discussing Szeg\H{o}'s contributions, one should mention that he was chair of the Stanford math department for over 20 years, turning the department from a provincial one to a world class one.  And his short term as chair of Washington University in St.~Louis left a long standing tradition that make it now, almost 90 years after he left, a serious center in harmonic and complex analysis.\\[-0.2cm]

\noindent Finally, I close with something I learned about Szeg\H{o} from Peter Lax (who sometimes said that Szeg\H{o} was a relative).  I have long been struck by the fact that classical analysts have not been sufficiently recognized by the U.S.~National Academy of Sciences and in particular that none of Aronszajn, Hartman or Wintner was elected.  After I wrote my book on Loewner's theorem, I gave a talk on the Wigner--von Neumann proof of that theorem at Courant and remarked that I was surprised that Loewner was never elected either.  Peter Lax said from the audience ``you know, neither was Szeg\H{o}''.  I was incredulous but checked overnight and the next day saw Peter at Courant tea and he told me:~``you know in those days, the NAS membership in math was restricted to only a few departments and Paul Cohen, who was elected after his Fields medal, was the first person in math elected at Stanford. He sought to get P\'{o}lya and Szeg\H{o} elected.  One day, I picked up the phone and heard \emph{Peter, this is Saunders, who is this P\'{o}lya}?'' (Saunders refers to Saunders MacLane)  P\'{o}lya was elected but Lax explained, by that time, Szeg\H{o} had Parkinson's disease which, given the attitude at the time, rendered him no longer a serious candidate.  And Peter completed what he said with ``and you know I think Szeg\H{o} was the deeper mathematician than P\'{o}lya'', a sentiment with which I agree.

\bibliographystyle{plain}
\bibliography{refbib}

\end{document}